\pdfoutput=1 
\RequirePackage[l2tabu, orthodox]{nag} 
\documentclass[a4paper,reqno]{article} 
\usepackage[stretch=10]{microtype} 
\usepackage[german,french,UKenglish]{babel}
\usepackage[T1]{fontenc} 
\usepackage[utf8]{inputenc} 
\usepackage{csquotes} 
\usepackage[backend=biber,style=alphabetic,giveninits=true, 
isbn=false,sorting=nyt,maxnames=4,doi=false,url=false]{biblatex}
\usepackage{xpatch} 
\usepackage{letltxmacro} 
\usepackage{perpage} 
\usepackage{soul} 
\usepackage{xparse} 
\usepackage{amsmath}
\usepackage{amsfonts}
\usepackage{amssymb}
\usepackage{amsthm}
\usepackage{mathrsfs}
\usepackage{mathdots}
\usepackage{enumitem}
\usepackage{calc}
\usepackage{tocloft} 
\usepackage{titlesec} 
\usepackage{textcomp} 
\usepackage{tikz}
\usepackage{tikz-cd} 
\usetikzlibrary{calc}
\usetikzlibrary{matrix,arrows}
\usetikzlibrary{decorations.markings}

\usepackage{verbatim} 
\usepackage{aliascnt} 
\usepackage[colorlinks=true,urlcolor=blue!70!black,citecolor=blue!60!black,linkcolor=blue!60!black,pdfauthor={R. van Dobben de Bruyn},pdftitle={Automorphisms of categories of schemes}]{hyperref} 

\titlespacing\section{0pt}{12pt plus 4pt minus 2pt}{0pt plus 2pt minus 2pt}
\titlespacing\subsection{0pt}{12pt plus 4pt minus 2pt}{0pt plus 2pt minus 2pt}
\titlespacing\subsubsection{0pt}{0pt plus 2pt minus 2pt}{0pt plus 2pt minus 2pt}

\titleformat{\section}[block]{\Large\bfseries\scshape\filcenter}{\thesection.}{1ex}{}
\titleformat{\subsection}{\large\scshape\filcenter}{\thesubsection}{1ex}{}

\makeatletter

\ExplSyntaxOn
\cs_new:Npn \white_text:n #1
  {
    \llap{\textcolor{white}{\the\SOUL@syllable}\hspace{\fp_eval:n{0.01*#1} em}}
    \llap{\textcolor{white}{\the\SOUL@syllable}\hspace{-\fp_eval:n{0.01*#1} em}}
  }
\NewDocumentCommand{\whiten}{ m }
    {
      \int_step_function:nnnN {1}{1}{#1} \white_text:n
    }
\ExplSyntaxOff

\NewDocumentCommand{ \varul }{ D<>{5} O{0.2ex} O{0.1ex} +m } {%
\begingroup
\setul{#2}{#3}%
\def\SOUL@uleverysyllable{%
   \setbox0=\hbox{\the\SOUL@syllable}%
   \ifdim\dp0>\z@
      \SOUL@ulunderline{\phantom{\the\SOUL@syllable}}%
      \whiten{#1}%
      \llap{%
        \the\SOUL@syllable
        \SOUL@setkern\SOUL@charkern
      }%
   \else
       \SOUL@ulunderline{%
         \the\SOUL@syllable
         \SOUL@setkern\SOUL@charkern
       }%
   \fi}%
    \ul{#4}%
\endgroup
}

\makeatother

\makeatletter
\tikzcdset{
  open/.code     = {\tikzcdset{hook, circled};},
  closed/.code   = {\tikzcdset{hook, slashed};},
  open'/.code    = {\tikzcdset{hook', circled};},
  closed'/.code  = {\tikzcdset{hook', slashed};},
  circled/.code  = {\tikzcdset{markwith = {\draw (0,0) circle (.375ex);}};},
  slashed/.code  = {\tikzcdset{markwith = {\draw[-] (-.4ex,-.4ex) -- (.4ex,.4ex);}};},
  markwith/.code ={
    \pgfutil@ifundefined%
    {tikz@library@decorations.markings@loaded}%
    {\pgfutil@packageerror{tikz-cd}{You need to say %
      \string\usetikzlibrary{decorations.markings} to use arrows with markings}{}}{}%
    \pgfkeysalso{/tikz/postaction = {
      /tikz/decorate,
      /tikz/decoration={markings, mark = at position 0.5 with {#1}}}
    }
  },
}
\makeatother

\newtheoremstyle{thms}{1em}{0pt}{\itshape}{}{\bfseries}{.}{ }{}
\theoremstyle{thms}
\newtheorem{Thm}{Theorem}[section]				

\newaliascnt{Prop}{Thm}							
\newtheorem{Prop}[Prop]{Proposition}
\aliascntresetthe{Prop}

\newaliascnt{Lemma}{Thm}							
\newtheorem{Lemma}[Lemma]{Lemma}
\aliascntresetthe{Lemma}

\newaliascnt{Cor}{Thm}						
\newtheorem{Cor}[Cor]{Corollary}
\aliascntresetthe{Cor}

\newaliascnt{Conj}{Thm}							

\aliascntresetthe{Conj}

\newaliascnt{Question}{Thm}						
\newtheorem{Question}[Question]{Question}
\aliascntresetthe{Question}

\newtheoremstyle{defs}{1em}{0pt}{}{}{\bfseries}{.}{ }{}
\theoremstyle{defs}
\newaliascnt{Rmk}{Thm}							
\newtheorem{Rmk}[Rmk]{Remark}
\aliascntresetthe{Rmk}

\newaliascnt{Fact}{Thm}							

\aliascntresetthe{Fact}

\newaliascnt{Def}{Thm}							
\newtheorem{Def}[Def]{Definition}
\aliascntresetthe{Def}

\newaliascnt{Ex}{Thm}								
\newtheorem{Ex}[Ex]{Example}
\aliascntresetthe{Ex}

\newaliascnt{Con}{Thm}							

\aliascntresetthe{Con}

\newaliascnt{Not}{Thm}							

\aliascntresetthe{Not}

\newaliascnt{Setup}{Thm}							

\aliascntresetthe{Setup}

\newaliascnt{Picture}{Thm}						

\aliascntresetthe{Picture}

\theoremstyle{thms}
\newtheorem{thm}{Theorem}
\newtheorem*{thm*}{Theorem}

\newtheorem*{lemma*}{Lemma}

\newtheorem*{cor*}{Corollary}
\newtheorem{cor}{Corollary}

\LetLtxMacro\oldproof\proof						
\renewcommand{\proof}[1][Proof]{\oldproof[#1]\unskip} 

\newcommand{\boldref}[1]{\autoref{#1}}

\newenvironment{itemize*} 
  {\begin{itemize}
    \setlength{\itemsep}{1em}
    \setlength{\parskip}{-1em}
    \setlength{\topsep}{0pt}
    \setlength{\partopsep}{0pt}}
  {\end{itemize}}
  
\newenvironment{enumerate*}
  {\begin{enumerate}
    \setlength{\itemsep}{1em}
    \setlength{\parskip}{-1em}
    \setlength{\topsep}{0pt}
    \setlength{\partopsep}{0pt}}
  {\end{enumerate}}
  
\setlist{itemsep=0em,topsep=0cm,partopsep=0em,parsep=\lineskip}

\setlist[enumerate]{label=\normalfont(\arabic*)}

\setlist[itemize]{leftmargin=1.3em}

\newcommand{\cat}[1]{\operatorname{\textbf{\textup{#1}}}} 
\newcommand{\Hom}{\operatorname{Hom}}
\newcommand{\End}{\operatorname{End}}

\newcommand{\Aut}{\operatorname{Aut}}
\newcommand{\Isom}{\operatorname{Isom}}
\newcommand{\ISOM}{\operatorname{\textbf{\textup{Isom}}}}
\newcommand{\Mor}{\operatorname{Mor}}

\newcommand{\Spec}{\operatorname{Spec}}
\newcommand{\SPEC}{\operatorname{\textbf{\textup{Spec}}}}
\newcommand{\Sch}{\cat{Sch}}

\newcommand{\NIL}{\operatorname{\textbf{\textup{Nil}}}}
\newcommand{\FEt}{\cat{F\'Et}}
\newcommand{\LRS}{\cat{LRS}}
\newcommand{\Ring}{\cat{Ring}}

\newcommand{\Set}{\cat{Set}}
\newcommand{\Top}{\cat{Top}}

\newcommand{\Qcoh}{\cat{Qcoh}}
\newcommand{\QCOH}{\cat{Qcoh}}
\newcommand{\Coh}{\cat{Coh}}
\newcommand{\Alg}{\cat{Alg}}

\newcommand{\Mod}{\cat{Mod}}

\newcommand{\Frac}{\operatorname{Frac}}

\newcommand{\coAbsep}{\cat{coAb}^{\operatorname{sep}}}

\newcommand{\op}{^{\operatorname{op}}}
\newcommand{\Sym}{\operatorname{Sym}}
\newcommand{\Pic}{\operatorname{Pic}}

\newcommand{\id}{\operatorname{id}}
\newcommand{\im}{\operatorname{im}}

\newcommand{\Ra}{\Rightarrow}

\newcommand{\LRa}{\Leftrightarrow}

\newcommand{\rsa}{\rightsquigarrow}

\newcommand{\punct}[1]{\makebox[0pt][l]{\,#1}} 

\newcommand{\N}{\mathbf N}
\newcommand{\Z}{\mathbf Z}
\newcommand{\Q}{\mathbf Q}
\newcommand{\R}{\mathbf R}
\newcommand{\C}{\mathbf C}
\newcommand{\F}{\mathbf F}
\renewcommand{\P}{\mathbf P}
\newcommand{\A}{\mathbf A}

\newcommand{\x}{^\times}

\DefineBibliographyStrings{UKenglish}{edition = {edition}} 
\AtEveryBibitem{\ifentrytype{book}{\clearfield{url}\clearfield{pages}}{}}
\renewbibmacro{in:}{}
\DeclareFieldFormat[book,incollection]{number}{\mkbibbold{#1}}
\DeclareFieldFormat[article]{volume}{\mkbibbold{#1}}
\DeclareFieldFormat*{title}{\mkbibitalic{#1}} 
\DeclareFieldFormat*{booktitle}{#1} 
\DeclareFieldFormat*{journaltitle}{#1}
\DeclareFieldFormat{postnote}{#1}
\DeclareFieldFormat{pages}{p. #1}
\DeclareFieldFormat[thesis]{type}{\bibstring{#1}\addcomma\addspace}
\DefineBibliographyStrings{english}{bibliography = {References},}
\AtEveryBibitem{\clearlist{publisher}\clearname{editor}\clearfield{location}}
\renewcommand*{\newunitpunct}{\addcomma\space}
\renewcommand*{\newblockpunct}{\addperiod\space}
\renewcommand*{\bibfont}{\footnotesize}

\renewcommand*{\mkbibnamefamily}[1]{\textsc{#1}} 
\renewcommand*{\mkbibnameprefix}[1]{\textsc{#1}}

\newbibmacro{titlelink}[1]{%
	\iffieldundef{doi}{%
		\iffieldundef{url}{%
			#1
		}{%
			\href{\thefield{url}}{#1}%
		}%
	}{%
		\href{https://dx.doi.org/\thefield{doi}}{#1}%
	}%
}
\DeclareFieldFormat{title}{\usebibmacro{titlelink}{\mkbibemph{#1}}}

\xpatchbibdriver{article}
  {\newunit
   \usebibmacro{note+pages}}
  {}
  {}{}
\renewbibmacro*{journal+issuetitle}{%
  \usebibmacro{journal}%
  \setunit*{\addspace}%
  \iffieldundef{series}
    {}
    {\newunit
     \printfield{series}%
     \setunit{\addspace}}%
  \usebibmacro{volume+number+eid}%
  \newunit
  \usebibmacro{note+pages}%
  \setunit{\addspace}%
  \usebibmacro{issue+date}%
  \setunit{\addcolon\space}%
  \usebibmacro{issue}%
  \newunit}
\xpatchbibdriver{inbook}
  {\newunit\newblock
  \usebibmacro{chapter+pages}}%
  {}%
  {}{}
\xpatchbibdriver{incollection}
  {\newunit\newblock
  \usebibmacro{chapter+pages}}%
  {}%
  {}{}
\xpatchbibdriver{inproceedings}
  {\newunit\newblock
  \usebibmacro{chapter+pages}}%
  {}%
  {}{}
\xapptobibmacro{maintitle+booktitle}
  {\printfield{pages}%
  \newunit}
  {}{}

\xpatchbibdriver{article}
  {\newunit\newblock
  \usebibmacro{in:}}%
  {\setunit{\addperiod\space}
  \usebibmacro{in:}}%
  {}{}
\xpatchbibdriver{book}
  {\newunit\newblock
  \usebibmacro{series+number}%
  \newunit\newblock
  \printfield{note}}%
  {\setunit{\addperiod\space}%
  \usebibmacro{series+number}
  \setunit{\addperiod\space}%
  \printfield{note}}%
  {}{}
\xpatchbibdriver{inbook}
  {\newunit\newblock
  \usebibmacro{maintitle+booktitle}}%
  {\setunit{\addperiod\space}
  \usebibmacro{maintitle+booktitle}}%
  {}{}
\xpatchbibdriver{incollection}
  {\newunit\newblock
  \usebibmacro{in:}}%
  {\setunit{\addperiod\space}}
  {}{}
\xpatchbibdriver{inproceedings}
  {\newunit\newblock
  \usebibmacro{in:}}%
  {\setunit{\addperiod\space}}%
  {}{}
\xpatchbibdriver{incollection}
  {\newunit\newblock
  \usebibmacro{series+number}}%
  {\setunit{\addperiod\space}
  \usebibmacro{series+number}}%
  {}{}
\xpatchbibdriver{thesis}
  {\newunit\newblock
  \printfield{type}}%
  {\setunit{\addperiod\space}
  \printfield{type}}%
  {}{}
\xpatchbibdriver{misc}
  {\newunit\newblock
  \printfield{howpublished}}%
  {\setunit{\addperiod\space}\newblock
  \printfield{howpublished}}%
  {}{}
\xpatchbibmacro{doi+eprint+url}
  {\newunit\newblock
  \iftoggle{bbx:eprint}}%
  {\setunit{\addperiod\space}
  \iftoggle{bbx:eprint}}%
  {}{}
\xpatchbibmacro{series+number}
  {\printfield{series}}
  {\iffieldundef{shortseries}
    {\printfield{series}}
    {\printfield{shortseries}}}
  {}{}

\begin{document}

\renewcommand{\sectionautorefname}{Section}
\renewcommand{\subsectionautorefname}{Subsection}		

\begin{center}
\vspace*{-2em}
\noindent\makebox[\linewidth]{\rule{14cm}{0.4pt}}
\vspace{.5em}

\noindent\makebox[\linewidth]{\LARGE{\textsc{\textbf{Automorphisms of categories of schemes}}}}

\vspace{1.2em}

{\large{\textsc{Remy van Dobben de Bruyn}}}

\vspace{.5em}
\rule{8cm}{0.4pt}
\vspace{3em}
\end{center}

\renewcommand{\abstractname}{\small\bfseries\scshape Abstract}

\begin{abstract}\noindent
Given two schemes $S$ and $S'$, we prove that every equivalence between $\Sch_S$ and $\Sch_{S'}$ comes from a unique isomorphism between $S$ and $S'$. This eliminates all Noetherian and finite type hypotheses from a result of Mochizuki \cite{Moch} and fully answers a programme set out by Brandenburg in a series of questions on MathOverflow in 2011 \cite{Brand1,Brand2,Brand3,Brand4}.
\end{abstract}


\phantomsection
\section*{Introduction}\label{Sec intro}
Let $\mathscr C$ and $\mathscr D$ be categories, and write $\ISOM(\mathscr C,\mathscr D)$ for the category whose objects are equivalences $F \colon \mathscr C \to \mathscr D$ and whose morphisms $\eta \colon F \to G$ are natural isomorphisms. Let $\Isom(\mathscr C,\mathscr D) = \pi_0(\ISOM(\mathscr C,\mathscr D))$ be its set of isomorphism classes, i.e.\ the set of equivalences $F \colon \mathscr C \to \mathscr D$ up to natural isomorphism. 
We will study $\ISOM(\mathscr C,\mathscr D)$ when $\mathscr C$ and $\mathscr D$ are categories of schemes. The main result is the following.

\begin{thm}\label{Thm main}
Let $S$ and $S'$ be schemes. Then the natural functor
\[
\Isom(S,S') \to \ISOM(\Sch_{S'},\Sch_S)
\]
is an equivalence (where $\Isom(S,S')$ is a discrete category).
\end{thm}

A similar result for the category of locally Noetherian schemes with finite type morphisms was proven by Mochizuki \cite[Thm.~1.7]{Moch}. Our result completely eliminates all Noetherian and finite type hypotheses. Because we do not have access to the same finite type techniques, the proof is almost entirely independent from Mochizuki's.

A special case of interest is the case $S = S' = \Spec \Z$. This gives a positive answer to a question by Brandenburg \cite{Brand1}:

\begin{cor}\label{Cor MO}
Let $F \colon \Sch \to \Sch$ be an equivalence. Then $F$ is isomorphic to the identity functor.
\end{cor}

While this paper was in preparation, \boldref{Cor MO} was obtained independently by Pohl \cite{Pohl}.

We also get a version for commutative rings:

\begin{thm}\label{Thm rings}
Let $R$ and $R'$ be commutative rings. Then the natural functor
\[
\Isom(R,R') \to \ISOM(\Alg_R,\Alg_{R'})
\]
is an equivalence.
\end{thm}

A first statement of this type was proven by Clark and Bergman in the case where $R = R'$ is a (commutative, unital) integral domain \cite[Thm.~5.5]{CB}. Their result is only on the level of $\pi_0$ (equivalences up to natural isomorphism), and does not address whether the $\ISOM$ category is a setoid. It also deals with categories of algebras that are not necessarily commutative or unital, under the same assumptions on $R = R'$. In the case of non-commutative algebras, one also gets the $(-)\op$ autoequivalence [\emph{loc.\ cit.}].

We imagine that \boldref{Thm rings} may have been known to experts, although we do not know a reference. The analogue of \boldref{Cor MO} for commutative (unital) rings follows from \cite{CB}; see also the answers to the post \cite{Bel} for a number of easy alternative proofs.

\subsection*{Strategy of proof}
The strategy of the proof of \boldref{Thm main} is to characterise certain properties of schemes and morphisms of schemes by purely categorical means. The general framework is that of \emph{categorical reconstruction}, in the sense of \boldref{Def reconstructed} below. Although this technique has been used before in many different contexts, as far as we know there is no systematic treatment in the literature. In \boldref{Sec categorical}, we recall the definitions and prove some basic lemmas, in particular in the setting of slice categories.

We then specialise to slice categories of schemes. The underlying set of a scheme is easily found as the isomorphism classes of simple subobjects (\boldref{Sec set}). To find the topology (\boldref{Sec topology}), we first relate locally closed immersions to regular monomorphisms (\boldref{Lem immersion}). Then we characterise spectra of valuation rings (\boldref{Prop valuation}), which with some work recovers the topology (\boldref{Prop topology}).

On the other hand, a variant of a standard argument due to Beck \cite[Ex.\ 8]{BeckThesis} recovers the category $\Qcoh(\mathcal O_S)$ as cogroup objects in the category $S/\Sch_S$ of $S$-schemes with a section (\boldref{Sec quasi-coherent}). The centre of $\Qcoh(\mathcal O_S)$ is $\Gamma(S,\mathcal O_S)$, and a sheafy version of this statement recovers the structure sheaf $\mathcal O_S$ (\boldref{Sec structure sheaf}). This finishes the reconstruction of a scheme isomorphic to $S$ from the slice category $\Sch_S$. By \boldref{Lem categorical reconstruction equivalence on ISOM}, a sufficiently functorial version of this immediately implies the main theorem.

The first proof we found of \boldref{Thm main} required an additional step in between \boldref{Sec topology} and \boldref{Sec quasi-coherent}, namely the reconstruction of affine morphisms. Pohl's argument for \boldref{Cor MO} shows that this step is not needed, because cogroups in $S/\Sch_S$ are nilpotent thickenings, hence automatically affine. The now omitted characterisation of affine morphisms gives some results of independent interest, which will appear separately.

\subsection*{Acknowledgements}
{\small
I am grateful to Gregor Pohl for the shortcut provided by his argument for \boldref{Cor MO}, bypassing entirely the characterisation of affine morphisms. I thank Raymond Cheng for many helpful discussions, and Bhargav Bhatt for helpful suggestions and further questions. Finally, I want to thank Johan de Jong for his encouragement to eliminate all remaining Noetherian and finite type hypotheses, as well as for providing most of \boldref{Lem closed immersion}.
}

\subsection*{Notation and conventions}
Categories will be written bold, like the categories of sets ($\Set$), topological spaces ($\Top$), rings ($\Ring$), $R$-algebras ($\Alg_R$), schemes ($\Sch$), locally ringed spaces ($\LRS$), etcetera. All rings and algebras will be commutative and unital. Given a morphism of schemes $f \colon X \to Y$, write $f^\#$ for the morphism of sheaves $\mathcal O_Y \to f_*\mathcal O_X$ on $Y$.

Given an object $X$ in a category $\mathscr C$, we write $\mathscr C/X$ for the slice category of objects over $X$, and $X/\mathscr C$ for the coslice category of objects under $X$. In the former case, we sometimes write $\mathscr C_{/X}$ or simply $\mathscr C_X$ if this causes no confusion. For example, $\Sch_S$ is the category $\Sch/S$ of schemes over a base scheme $S$. Note however that $\Alg_R$ is by convention the \emph{coslice} category $R/\Ring$.

To deal with set-theoretic issues arising in the formation of $\ISOM(\Sch_{S'}, \Sch_S)$ and other category theoretic constructions, one should either work with universes, or use `sufficiently large' small categories of schemes, cf.\ e.g.\ \cite[Tag \href{https://stacks.math.columbia.edu/tag/000J}{000J}]{Stacks}.

\numberwithin{equation}{section}

\section{Functorial reconstruction}\label{Sec categorical}
The strategy of the proof of \boldref{Thm main} is to show that certain properties of schemes are \emph{categorical}, in the sense of \boldref{Def categorical} below.

The mathematical foundations needed to define categorical reconstruction lie in the intersection of model theory and category theory. The notion appears to originate in \cite[\S2]{CB}, but we don't know a systematic study in the literature, and it is not always possible to give detailed references. We will illustrate the definitions with examples from algebra and topology.

\begin{Def}\label{Def categorical}
Let $\mathscr C$ be a category. Then a property $\mathcal P$ of objects $X \in \mathscr C$ or of morphisms $f \colon X \to Y$ in $\mathcal C$ is \emph{categorical} if it is definable in terms of morphisms, compositions, and equality of morphisms. This in particular implies that if $F \colon \mathscr C \to \mathscr D$ is an equivalence, then $\mathcal P(X) \LRa \mathcal P(FX)$ for an object $X$ in $\mathscr C$ (resp.~$\mathcal P(f) \LRa \mathcal P(Ff)$ for a morphism $f$ in $\mathscr C$).
\end{Def}

\begin{Rmk}
However, our language does not have a predicate for \emph{equality} (rather than isomorphism) of objects. For example, the property that an object is the only one in its isomorphism class should not be a categorical one, as it is not stable under equivalence of categories.
\end{Rmk}

\begin{Def}\label{Def reconstructed}
Let $F \colon \mathscr C \to \mathscr D$ be a functor to a concretely definable category $\mathscr D$. Then a \emph{categorical reconstruction} of $F$ is a pair $(F', \eta)$ consisting of a functor $F' \colon \mathscr C \to \mathscr D$ that is definable in terms of morphisms of $\mathscr C$ satisfying categorical properties, together with a natural isomorphism $\eta \colon F \to F'$.
\end{Def}

Here, a \emph{concretely definable} category $\mathscr D$ (used loosely) means a category of some collection of sets satisfying some relations\footnote{More precisely, there should be a predicate that takes a set and says whether it's in $\mathscr D$; see e.g.\ \cite[Tag \href{https://stacks.math.columbia.edu/tag/0009}{0009}]{Stacks} for the case $\mathscr D = \Sch$.}. For example, $\Sch$ is concretely definable, because a scheme is a topological space with a sheaf of rings on it. 
%
However, a category abstractly equivalent to $\Sch$ need not be concretely definable.

Note that the roles of $\mathscr C$ and $\mathscr D$ are rather asymmetric in the definition of categorical reconstruction. The classical case is $\mathscr D = \Set$:

\begin{Ex}
The forgetful functor $\Top \to \Set$ can be reconstructed categorically. Indeed, if $X$ is a topological space, then points of $X$ are in bijection with $\Mor(*,X)$, where $*$ is any terminal object.
\end{Ex}

\begin{Lemma}\label{Lem auto-equivalence}
Let $F \colon \mathscr C \to \mathscr D$ be a functor that can be reconstructed categorically. If $G \colon \mathscr C \to \mathscr C$ is an auto-equivalence, then $FG \cong F$.
\end{Lemma}

\begin{proof}
For any $X \in \mathscr C$, the auto-equivalence $G$ takes the data used to define $F(X)$ to the data used to define $F(G(X))$.
\end{proof}
\vskip-\lastskip
It turns out to be very powerful to reconstruct the identity functor $\mathscr C \to \mathscr C$ of a concretely definable category:

\begin{Cor}
If the identity functor $\id \colon \mathscr C \to \mathscr C$ of a concretely definable category can be reconstructed categorically, then $\Aut(\mathscr C) = 1$.\hfill\qed
\end{Cor}

For example, by the following lemma we get $\Aut(\Top) = 1$.

\begin{Lemma}
The identity $\Top \to \Top$ can be reconstructed categorically.
\end{Lemma}

\begin{proof}
The property that $X \in \Top$ is (isomorphic to) the Sierpi\'nski space $S$ is categorical: it is the unique two-point space for which the swap is not continuous, i.e.~$\#\Mor(*,S) = 2$ and $\Aut(S) = 1$.

Then the functor $\Top\op \to \Set$ given by $(X,\mathcal T) \mapsto \mathcal T$ can also be reconstructed categorically, as it is corepresented by $S$. The open point $\eta \in S$ is characterised as the unique point such that the inclusion
\begin{align*}
\iota_X \colon \Mor(X,S) &\hookrightarrow \mathcal P(\Mor(*,X))\\
f &\mapsto \{g \in \Mor(*,X)\ |\ fg = \eta\}
\end{align*}
endows $\Mor(*,X)$ with a topology for all $X$ (equivalently, the system of sets $\iota_X(f)$ for $f \in \Mor(X,S)$ is closed under arbitrary unions).

The functorial association
\[
X \mapsto \Big(\Mor(*,X),\iota_X(\Mor(X,S))\Big)
\]
recovers a topological space naturally homeomorphic to $X$.
\end{proof}
\vskip-\lastskip
The categorical reconstruction of the forgetful functor $\Sch_S \to \Sch$ follows a similar strategy. We start by reconstructing the forgetful functor $\Sch_S \to \Set$ (\boldref{Cor basic properties}), then upgrade this to $\Sch_S \to \Top$ (\boldref{Prop topology}), and finally we reconstruct $\Sch_S \to \Sch$ (\boldref{Thm forgetful functor}).

There are many examples of categorical reconstruction theorems in the literature; a sample of well-known results includes the Neukirch--Uchida theorem \cite{Neu1,Neu2,Uch1,Uch2}, the Gabriel--Rosenberg theorem \cite{Gabriel,Ros,BrandRosenberg}, the Bondal--Orlov theorem \cite{BonOrl}, and Mochizuki's results in anabelian geometry \cite{MochAnab1}, \cite{MochAnab2} and for Noetherian schemes and log schemes \cite{Moch}. The most general setup is as follows.

\begin{Question}
Let $\mathscr C$ be a finitely complete category, and let $\mathscr S \to \mathscr C$ be a fibred category. If $\mathscr S_X$ and $\mathscr S_Y$ are equivalent for $X, Y \in \mathscr C$, is it true that $X \cong Y$?
\end{Question}

\begin{Rmk}
In general it is too much to expect that the natural map
\[
\Isom(X,Y) \to \Isom(\mathscr S_Y,\mathscr S_X)
\]
is an isomorphism. For example, although Gabriel \cite{Gabriel} (resp.~Rosenberg \cite{Ros}, \cite{BrandRosenberg}) reconstruct a Noetherian scheme (resp.~quasi-separated scheme) $X$ from its category $\Qcoh(\mathcal O_X)$, the latter can have extra endomorphisms not coming from $X$: if $\mathscr L$ is a nontrivial line bundle, then $-\otimes \mathscr L$ is an auto-equivalence of $\Qcoh(\mathcal O_X)$ that does not come from an automorphism of $X$ (e.g. since it does not fix $\mathcal O_X$).

Similarly, Bondal and Orlov \cite{BonOrl} reconstruct a smooth projective variety $X$ with ample canonical or anti-canonical bundle from its derived category $D(X) = D^b_{\Coh}(X)$, but the derived category has a shift $[1] \colon D(X) \to D(X)$ that is not induced by an automorphism of $X$. In fact, they prove \cite[Thm.\ 3.1]{BonOrl} that $\Aut(D(X)) \cong \Aut(X) \ltimes (\Pic(X) \oplus \Z)$, where $\Z$ acts by shifting.
\end{Rmk}

However, for slice categories $\mathscr C_{/X}$ we have the following result, which is often implicitly reproved in applications. For example, this argument applies to the Noetherian version \cite{Moch}, as well as some anabelian situations (taking slice categories $\FEt_{/X}$ in the category of finite \'etale morphisms of schemes).

\begin{Lemma}\label{Lem categorical reconstruction equivalence on ISOM}
Let $\mathscr C$ be a finitely complete category. Assume that the forgetful functors $F_X \colon \mathscr C_{/X} \to \mathscr C$ can be reconstructed categorically from $\mathscr C_{/X}$, by a formula that does not depend on $X$. Then the natural functors
\begin{align*}
\Psi \colon \Isom(X,Y) &\to \ISOM(\mathscr C_{/Y},\mathscr C_{/X})\\
f &\mapsto f^*
\end{align*}
for $X, Y \in \mathscr C$ are equivalences of categories (where the left hand side is viewed as a discrete category).
\end{Lemma}

That is, $\ISOM(\mathscr C_{/Y},\mathscr C_{/X})$ is a setoid whose isomorphism classes are $\Isom(X,Y)$.

\begin{proof}
For simplicity, assume the categorical reconstructions of the forgetful functors are \emph{equal} (rather than isomorphic) to $F_X$ and $F_Y$; this does not affect the argument. If $F \colon \mathscr C_{/Y} \to \mathscr C_{/X}$ is an equivalence, then there exists a natural isomorphism $\eta \colon F_XF \stackrel\sim\to F_Y$. For any morphism $B \to Y$, let its image under $F$ be $A \to X$. Then naturality of $\eta$ gives a commutative diagram
\begin{equation}
\begin{tikzcd}\label{Dia essentially surjective}
A \ar{r}{\sim}[swap]{\eta_B}\ar{d} & B \ar{d}\\
X \ar{r}{\sim}[swap]{\eta_Y} & Y\punct{,}
\end{tikzcd}
\end{equation}
inducing functorial isomorphisms $\beta_B \colon A \stackrel\sim\to B \times_Y X$. Setting $f = \eta_Y \colon X \stackrel\sim\to Y$, we conclude that $\beta$ gives a natural isomorphism $\beta \colon F \stackrel\sim\to f^*$, showing that $\Psi$ is essentially surjective.

For fully faithfulness, it suffices to show that $\Aut(\id_{\mathscr C_{/X}}) = 1$ for all $X \in \mathscr C$. Let $\alpha \colon \id_{\mathscr C_{/X}} \stackrel\sim\to \id_{\mathscr C_{/X}}$ be an automorphism. For every morphism $f \colon A \to B$ in $\mathscr C_{/X}$, naturality of $\alpha$ gives a commutative diagram
\[
\begin{tikzcd}
A \ar{d}[swap]{f}\ar{r}{\sim}[swap]{\alpha_A} & A\ar{d}{f}\\
B \ar{d}\ar{r}{\sim}[swap]{\alpha_B} & B\ar{d}\\
X \ar{r}{\sim}[swap]{\alpha_X} & X\punct{.}
\end{tikzcd}
\]
Applying this to all morphisms used in the categorical reconstruction of $F_X$, we conclude that $F_X(\alpha) = \id_{F_X}$. Since $F_X$ is faithful, this implies that $\alpha = \id$.
\end{proof}

\begin{Ex}
The forgetful functor $\Set_{/X} \to \Set$ can be reconstructed categorically: for a set $A$ over $X$, the points of $A$ correspond to (isomorphism classes of) simple subobjects of $A$. The proof of \boldref{Lem categorical reconstruction equivalence on ISOM} then reads as follows.
\begin{itemize}
\item If $F \colon \Set_{/Y} \to \Set_{/X}$ is an equivalence, then looking at the simple subobjects of the terminal objects gives an isomorphism $\eta_Y \colon X \to Y$, and diagram (\ref{Dia essentially surjective}) shows that $F \cong \eta_Y^*$.
\item If $\alpha \colon \id_{\Set_{/X}} \to \id_{\Set_{/X}}$ is an automorphism, then for any $A \in \Set_{/X}$ and any simple subobject $a \to A$, naturality of $\alpha$ gives the diagram
\[
\begin{tikzcd}
a \ar{d}\ar{r}{\sim}[swap]{\alpha_a} & a\ar{d}\\
A \ar{r}{\sim}[swap]{\alpha_A} & A\punct{,}
\end{tikzcd}
\]
showing that $\alpha_A$ fixes all points of $A$, i.e. $F_X(\alpha_A)$ agrees with $\id_{F_X(A)}$.
\end{itemize}
\end{Ex}

\begin{Rmk}
To prove \boldref{Thm main}, it therefore suffices to show that the forgetful functor $\Sch_S \to \Sch$ can be reconstructed categorically from $\Sch_S$.
\end{Rmk}

\section{Underlying set}\label{Sec set}
We begin by recovering the underlying set of an $S$-scheme $X$ from categorical information in $\Sch_S$; see \boldref{Cor basic properties} below.

\begin{Lemma}\label{Lem field}
The simple objects in $\Sch_S$ are the spectra of fields. \qed
\end{Lemma}

Here, a \emph{simple} objects is a nonempty object $X \in \Sch_S$ whose only subobjects are $\varnothing$ and $X$.

\begin{Lemma}\label{Lem points}
Let $X \in \Sch_S$. If $\Spec k \to X$ is a monomorphism from the spectrum of a field, then there exists a unique point $x \in X$ and a unique isomorphism $\Spec k \cong \Spec \kappa(x)$ over $X$. \qed
\end{Lemma}

\begin{Lemma}\label{Lem connected}
Let $X \in \Sch_S$. Then $X$ is connected if and only if $X$ cannot be written as a coproduct of two nonempty $S$-schemes. \qed
\end{Lemma}

To summarise the results of this section:

\begin{Cor}\label{Cor basic properties}
The following properties on objects $X \in \Sch_S$ are categorical:
\begin{enumerate}
\item $X$ is the spectrum of a field;\label{Cat field}
\item $X$ is connected.\label{Cat connected}
\end{enumerate}
Moreover, the forgetful functor $\Sch_S \to \Set$ can be reconstructed categorically.
\end{Cor}

\begin{proof}
Statements \ref{Cat field} and \ref{Cat connected} follow from \boldref{Lem field} and \boldref{Lem connected}. For the final statement, by \boldref{Lem points} we may take the functor that takes $X$ to the set of (isomorphism classes of) simple subobjects of $X$.
\end{proof}

\section{Topology}\label{Sec topology}
Having reconstructed the underlying set $|X|$ of an $S$-scheme $X$, we turn next to the topology on $|X|$. We first find locally closed immersions (\boldref{Lem immersion}), then describe spectra of valuation rings (\boldref{Prop valuation}), and finally recover the topology (\boldref{Prop topology}).

\begin{Lemma}\label{Lem immersion}
Let $f \colon X \to Y$ be a morphism in $\Sch_S$. Then $f$ is a (locally closed) immersion if and only if $f$ can be written as a composition of two regular monomorphisms. Thus, the property that $f$ is an immersion is categorical.
\end{Lemma}

\begin{proof}
If $f$ is the equaliser of $g_{1,2} \colon Y \rightrightarrows Z$, then the diagram
\begin{equation*}
\begin{tikzcd}[column sep=2.8em]
X \ar{d}[swap]{f}\ar{r} & Z \ar{d}{\Delta_Z}\\
Y\ar{r}{(g_1,g_2)} & Z \times Z
\end{tikzcd}
\end{equation*}
is a pullback. Hence, $f$ is an immersion since $\Delta_Z$ is \cite[Tag \href{https://stacks.math.columbia.edu/tag/01KJ}{01KJ}]{Stacks}. Thus, any regular monomorphism is an immersion, so the same goes for a composition of two regular monomorphisms \cite[Tag \href{https://stacks.math.columbia.edu/tag/02V0}{02V0}]{Stacks}.

Conversely, every immersion $f \colon X \to Y$ factors as $X \to U \to Y$, with $X \to U$ closed and $U \to Y$ open. Thus, it suffices to show that if $f$ is either an open immersion or a closed immersion, then $f$ is a regular monomorphism. 

If $f$ is an open immersion, then let $Z = Y \amalg_X Y$ be two copies of $Y$ glued along $X$, and let $g_{1,2} \colon Y \rightrightarrows Z$ be the two inclusions. Then $X$ is the equaliser of $Y \rightrightarrows Z$.

If $f$ is a closed immersion, then let $Z = Y \amalg_X Y$ (which exists and is for example described in \cite[Tag \href{https://stacks.math.columbia.edu/tag/0B7M}{0B7M}]{Stacks}). Once again, $X$ is the equaliser of $Y \rightrightarrows Z$. This proves the first statement, and the second is immediate.
\end{proof}

\begin{Rmk}
The proof actually shows that we may replace `composition of two regular monomorphisms' by `composition of two effective monomorphisms'.

One can combine the arguments in the open and closed cases to show that every immersion that factors as an open immersion inside a closed immersion is also a regular monomorphism \cite{Brand4} (but this no longer gives an effective monomorphism). As far as we know, the classification of regular monomorphisms in $\Sch$ is still open. For example, it does not appear to be known whether every regular monomorphism factors as open inside closed; see \cite{Brand4}.
\end{Rmk}


\begin{Lemma}\label{Lem reduced}
Let $X \in \Sch_S$. Then $X$ is reduced if and only if every immersion $Z \to X$ that induces a bijection on points $|Z| \to |X|$ is an isomorphism. Thus, the property that $X$ is reduced is categorical.
\end{Lemma}

\begin{proof}
An immersion whose underlying set map is closed is a closed immersion \cite[Tag \href{https://stacks.math.columbia.edu/tag/01IQ}{01IQ}]{Stacks}. This proves the first statement, and the second is immediate because immersions are categorical (\boldref{Lem immersion}) and the forgetful functor $\Sch_S \to \Set$ can be reconstructed categorically (\boldref{Cor basic properties}).
\end{proof}

\begin{Rmk}
A different characterisation was given by Moret-Bailly's partial answer \cite{MBAns} to \cite{Brand1}: an object $X \in \Sch_S$ is reduced if and only if the natural map $\coprod_{x \in X} \Spec \kappa(x) \to X$ is an epimorphism. Although this criterion is arguably more elementary, we find \boldref{Lem reduced} more intuitive.
\end{Rmk}

\begin{Prop}\label{Prop valuation}
Let $X \in \Sch_S$, and let $x \in X$ be a point. Then $(X,x)$ is isomorphic to $(\Spec R,\mathfrak m)$ for a valuation ring $R$ with maximal ideal $\mathfrak m$ if and only if all of the following hold:
\begin{enumerate}
\item $X$ is reduced and connected;\label{Item reduced and connected}
\item the category of immersions $Z \to X$ containing $x$ is a linear order;\label{Item totally ordered}
\item there exists a set $V \subseteq |X|$ that is the support of infinitely many immersions $Z \to X$ containing $x$.\label{Item support}
\end{enumerate}
In particular, the property that $(X,x) \cong (\Spec R,\mathfrak m)$ for a valuation ring $R$ with maximal ideal $\mathfrak m$ is categorical.
\end{Prop}

\begin{proof}
If $(X,x) \cong (\Spec R,\mathfrak m)$ with $R$ a valuation ring, then $X$ is reduced and connected since $R$ is a domain. Moreover, the only open subset $U \subseteq X$ containing $x$ is $X$, hence every immersion $Z \to X$ containing $x$ is closed.

Then the immersions containing $x$ are linearly ordered because the ideals of $R$ are linearly ordered. Finally, if $r \in R$ is a non-unit, then the ideals $(r) \supsetneq (r^2) \supsetneq \ldots$ all have the same underlying closed set $V \subseteq |X|$. Thus, \ref{Item reduced and connected}, \ref{Item totally ordered}, and \ref{Item support} hold if $(X,x) \cong (\Spec R,\mathfrak m)$ for a valuation ring $R$ with maximal ideal $\mathfrak m$.

Conversely, suppose $(X,x)$ satisfies \ref{Item reduced and connected}, \ref{Item totally ordered}, and \ref{Item support}. Firstly, we show that every point specialises to $x$. For a point $y \in X$, write $V(y)$ for the closure of $y$, and $D(y)$ for its complement.

First assume $y \in X$ neither specialises nor generalises to $x$, so $x \in D(y)$ and $y \in D(x)$. Define the locally closed sets $U = D(y)$ and $V = V(x) \cup V(y)$, both of which contain $x$. We have $y \in V$ and $y \not\in U$, hence $V \nsubseteq U$, so $U \subseteq V$ by \ref{Item totally ordered}. This means that $D(y) \subseteq V(x)$, hence $D(x) \cap D(y) = \varnothing$.

On the other hand, if $U = D(x) \triangle D(y)$ (the symmetric difference) and $V = V(x)$, then $U$ and $V$ both contain $x$, but $y$ is in $U$ and not in $V$. Hence, $U \nsubseteq V$, so \ref{Item totally ordered} implies $V \subseteq U$. This means that $V(x) \subseteq D(y)$, hence $D(x) \cup D(y) = X$. Thus, $X = D(x) \amalg D(y)$, contradicting connectedness of $X$.

Thus, every $y \in X$ either specialises to $x$ or generalises to $x$. Suppose that there exist strict specialisations $y \rsa x \rsa z$. Then the locally closed sets $U = D(z)$ and $V = V(y)$ contain $x$. But we have $y \in U$ but $y \not\in V$, and $z \in V$ but $z \not\in U$. This contradicts \ref{Item totally ordered}, and we conclude that one of the following holds:
\begin{itemize}
\item Every point $y \in X$ specialises to $x$;
\item Every point $y \in X$ generalises to $x$.
\end{itemize}
In the first case, the only open containing $x$ is $X$, so every immersion $Z \to X$ containing $x$ is closed. Similarly, in the second case every immersion $Z \to X$ containing $x$ is open. But since $X$ is reduced, the latter implies there is a unique scheme structure on every immersion $Z \to X$ containing $x$, contradicting \ref{Item support}. So we conclude that every point $y \in X$ specialises to $x$, and every immersion $Z \to X$ containing $x$ is closed.

Condition \ref{Item totally ordered} now also implies that for any $y,z \in X$, we have either $V(y) \subseteq V(z)$ or $V(z) \subseteq V(y)$, hence $z \rsa y$ or $y \rsa z$. In particular, $X$ has a unique generic point. Since $X$ is also assumed reduced, this implies that $X$ is integral. Since $X$ is the only open containing $x$, we conclude that $X$ is affine; say $X = \Spec R$ for some domain $R$. Finally, \ref{Item totally ordered} implies that the ideals in $R$ are linearly ordered, so $R$ is a valuation ring or a field. The maximal ideal $\mathfrak m \subseteq R$ corresponds to $x$, and $R$ is not a field by \ref{Item support}. This shows that \ref{Item reduced and connected}, \ref{Item totally ordered}, and \ref{Item support} imply that $(X,x) \cong (\Spec R, \mathfrak m)$ for a valuation ring $R$ with maximal ideal $\mathfrak m$.

The final statement follows as immersions are categorical (\boldref{Lem immersion}), reducedness and connectedness are categorical (\boldref{Lem reduced}, \boldref{Lem connected}), and the forgetful functor $\Sch_S \to \Set$ can be reconstructed categorically (\boldref{Cor basic properties}).
\end{proof}

\begin{Lemma}\label{Lem specialisation}
Let $X \in \Sch_S$, and let $x,y \in X$. Then $x \rsa y$ if and only if there exists a morphism $f \colon Z \to X$ where $Z$ is the spectrum of a valuation ring $R$ such that $f(\mathfrak m) = y$ and $x \in \im(f)$. In particular, the property that $x$ specialises to $y$ is categorical.
\end{Lemma}

\begin{proof}
Assume that such a morphism $f \colon Z \to X$ exists. Then $f^{-1}(V(x))$ is a nonempty closed subset of $Z$, hence contains $\mathfrak m$. This forces $y \in V(x)$, so $x \rsa y$.

Conversely, assume $x \rsa y$. The case $x = y$ is trivial, so we may assume $x \neq y$. Let $U \subseteq X$ be an affine open neighbourhood of $y$. Since $U$ is stable under generalisation, it contains $x$. If $U = \Spec A$ and $x,y \in U$ correspond to primes $\mathfrak p, \mathfrak q \subseteq A$ respectively, then we get a map $A \to B = (A/\mathfrak p)_{\mathfrak q}$. Since $\mathfrak p \subsetneq \mathfrak q$, the local domain $B$ is not a field. Hence, there exists a valuation ring $R \subseteq \Frac B$ dominating $B$ \cite[Tag \href{https://stacks.math.columbia.edu/tag/00IA}{00IA}]{Stacks}. Thus, the point $\mathfrak m_R \in \Spec R$ maps to $y \in X$, and $(0) \in \Spec R$ maps to $x \in X$.
%
This proves the first statement, and the final statement follows from \boldref{Cor basic properties} and \boldref{Prop valuation}.
\end{proof}
\vskip-\lastskip
It is tempting at this point to try to classify closed immersions as immersions (\boldref{Lem immersion}) that are closed under specialisation (\boldref{Lem specialisation}). However, this is not true in general; see \boldref{Ex closed immersion}. Instead, we use the following:

\begin{Lemma}\label{Lem closed immersion}
Given a map $f \colon X \to Y$ in $\Sch_S$, the following are equivalent:
\begin{enumerate}
\item $f$ is a closed immersion;\label{Item closed immersion}
\item $f$ is an immersion, and for every $Y' \to Y$ and every closed point $y \in Y'$ not in the image of the base change $f' \colon X' \to Y'$, the map $X' \amalg \Spec \kappa(y) \to Y'$ is an immersion.\label{Item disjoint union with closed point}
\end{enumerate}
In particular, closed immersions are categorical.
\end{Lemma}

\begin{proof}
If $f$ is a closed immersion, then clearly $f$ is an immersion, and for every $Y' \to Y$, the base change $f' \colon X' \to Y'$ is a closed immersion. If $y \not\in \im(f')$ is a closed point, then $g \colon X' \amalg \Spec \kappa(y) \to Y'$ is in fact a \emph{closed} immersion by the Chinese remainder theorem. This proves $\ref{Item closed immersion} \Ra \ref{Item disjoint union with closed point}$.

Conversely, assume \ref{Item disjoint union with closed point} holds. Since $f$ is an immersion, it suffices to show that the image of $f$ is closed \cite[Tag \href{https://stacks.math.columbia.edu/tag/01IQ}{01IQ}]{Stacks}. This can be done on a cover by affine opens, so let $Y' \subseteq Y$ be an affine open, with preimage $X' \subseteq X$.
%
If $y \in Y'$ is a closed point not in $X'$, then the map $X' \amalg \Spec \kappa(y) \to Y'$ is an immersion by \ref{Item disjoint union with closed point}. In particular, the subspace topology on $X' \cup \{y\}$ is the disjoint union of $X'$ and $\{y\}$, hence there exists an open $U \subseteq Y'$ such that $U \cap (X' \cup \{y\}) = \{y\}$. This means that $U$ is an open neighbourhood of $y$ with $U \cap X' = \varnothing$, so $y \not \in \overline{X'}$. Therefore the closed set $Z = \overline{X'}\setminus X' \subseteq Y'$ has no closed points, which implies $Z = \varnothing$ since $Y'$ is affine. Hence, $X'$ is closed in $Y'$, which proves $\ref{Item disjoint union with closed point} \Ra \ref{Item closed immersion}$.

The final statement follows since immersions are categorical (\boldref{Lem immersion}), the forgetful functor $\Sch_S \to \Set$ can be reconstructed categorically (\boldref{Cor basic properties}), and closed points are categorical (\boldref{Lem specialisation}).
\end{proof}

\begin{Prop}\label{Prop topology}
The forgetful functor $\Sch_S \to \Top$ can be reconstructed categorically.
\end{Prop}

\begin{proof}
By \boldref{Cor basic properties}, the forgetful functor $\Sch_S \to \Set$ can be reconstructed categorically. Moreover, a subset $U \subseteq |X|$ is open if and only if $U^c$ is the support of a closed immersion. The result now follows from \boldref{Lem closed immersion}.
\end{proof}

\begin{Cor}\label{Cor topological}
The following properties for a morphism $f \colon X \to Y$ in $\Sch_S$ are categorical:
\begin{enumerate}
\item $f$ is quasi-compact;\label{Cat qc}
\item $f$ is quasi-separated;\label{Cat qs}
\item $f$ is separated.\label{Cat separated}
\end{enumerate}
\end{Cor}

\begin{proof}
Statement 
\ref{Cat qc} is immediate from \boldref{Prop topology}. Statements \ref{Cat qs} and \ref{Cat separated} follow from \ref{Cat qc} and \boldref{Lem closed immersion} respectively, since $\Delta_{X/Y} \colon X \to X\times_Y X$ is functorially associated to $f$.
\end{proof}

\begin{Lemma}\label{Lem open immersion}
Let $f \colon X \to Y$ be a morphism. Then $f$ is an open immersion if and only if there exists an open set $U \subseteq |Y|$ such that $f$ is terminal among morphisms $Z \to Y$ landing in $U$. In particular, open immersions are categorical.
\end{Lemma}

\begin{proof}
The first statement is clear; the second follows from \boldref{Prop topology}.
\end{proof}
\vskip-\lastskip
If $X$ is Noetherian, then it suffices to replace the condition that $U \subseteq |Y|$ be open by the condition that $f$ is an immersion. This is not true for arbitrary schemes, as can be seen by the following example.

\begin{Ex}\label{Ex closed immersion}
Let $R$ be a ring admitting a pure ideal $I \subseteq R$ that is not finitely generated. For example, let $R = C^\infty(\R)$, and let $I$ be the functions vanishing in a neighbourhood of $0$ \cite[Tag \href{https://stacks.math.columbia.edu/tag/052H}{052H}]{Stacks}. Then there exists a multiplicative set $S \subseteq R$ such that $R/I \cong R[S^{-1}]$ as $R$-algebras \cite[Tag \href{https://stacks.math.columbia.edu/tag/04PS}{04PS}]{Stacks}.

Then the closed subset $V(I) \subseteq \Spec R$ is closed and stable under generalisation, but not open \cite[Tags \href{https://stacks.math.columbia.edu/tag/04PU}{04PU} and \href{https://stacks.math.columbia.edu/tag/05KK}{05KK}]{Stacks}. Similarly, the complementary open $D(I)$ is open and closed under specialisation, but not closed. Thus, we cannot get the closed immersions directly from \boldref{Lem immersion} and \boldref{Lem specialisation} as the immersions that are closed under specialisation.

Finally, the map $\Spec R/I \to \Spec R$ is terminal for maps landing in $V(I)$, because of the isomorphism $R/I \cong R[S^{-1}]$ and the universal property of rings of fractions. This shows that in \boldref{Lem open immersion}, we cannot replace the assumption that $U$ is open by the weaker assumption that $U$ is locally closed.
\end{Ex}

\section{Quasi-coherent sheaves}\label{Sec quasi-coherent}
We use a variant of the classical cogroup argument due to Beck \cite[Ex.\ 8]{BeckThesis} to recover the category of quasi-coherent sheaves on $X$ from the category $\Sch_X$.

\begin{Def}
For $X \in \Sch_S$, write $\Qcoh(\mathcal O_X)$ for the category of quasi-coherent $\mathcal O_X$-modules. Given an object $\mathscr F \in \Qcoh(\mathcal O_X)$, denote by $\Sym^{\leq 1} \mathscr F$ the quasi-coherent $\mathcal O_X$-algebra $\Sym^* \mathscr F/\Sym^2 \mathscr F \cong \mathcal O_X \oplus \mathscr F$, and by $\NIL_X(\mathscr F)$ the $X$-scheme $\SPEC_X(\Sym^{\leq 1} \mathscr F)$. It is the \emph{nilpotent thickening of $X$ by $\mathscr F$}.

The natural $\mathcal O_X$-algebra map $\Sym^{\leq 1} \mathscr F \twoheadrightarrow \mathcal O_X = \Sym^{\leq 0} \mathscr F$ is a retraction of $\mathcal O_X \to \Sym^{\leq 1} \mathscr F$, hence induces a section $\sigma \colon X \to \NIL_X(\mathscr F)$. This realises $\NIL_X(\mathscr F)$ as an object of the coslice category $X/\Sch_X$.
\end{Def}

\begin{Rmk}\label{Rmk pushout of schemes}
For $X \stackrel{\sigma_i}\to Y_i \to X$ in $X/\Sch_X$ (for $i \in \{1,2\}$), the pushout $Y_1 \amalg_X Y_2$ exists in $\Sch_X$ when each $\sigma_i$ is a closed immersion \cite[Tag \href{https://stacks.math.columbia.edu/tag/0B7M}{0B7M}]{Stacks}. Moreover, $|Y_1 \amalg_X Y_2|$ is the topological pushout $|Y_1| \amalg_{|X|} |Y_2|$, with the sheaf of rings given by $\mathcal O_{Y_1} \times_{\mathcal O_X} \mathcal O_{Y_2}$.

In particular, if $\mathscr F, \mathscr G \in \Qcoh(\mathcal O_X)$, then $\NIL_X(\mathscr F) \amalg_X \NIL_X(\mathscr G)$ exists, and the explicit description gives
\[
\NIL_X(\mathscr F) \underset X\amalg \NIL_X(\mathscr G) \cong \NIL_X(\mathscr F \oplus \mathscr G).
\]
Therefore, the functor $\NIL_X(-) \colon \Qcoh(\mathcal O_X)\op \to X/\Sch_X$ preserves finite coproducts. In particular, addition $\mathscr F \oplus \mathscr F \to \mathscr F$ induces a comultiplication
\[
c \colon \NIL_X(\mathscr F) \to \NIL_X(\mathscr F) \underset X\amalg \NIL_X(\mathscr F),
\]
and similarly for inversion $\iota$, making $(\NIL_X(\mathscr F),c,\iota)$ into an abelian cogroup object in $X/\Sch_X$. 
%
Here a cogroup object in $X/\Sch_X$ means for the monoidal structure given by the pushout $- \amalg_X -$; in particular if $Y$ is a cogroup in $X/\Sch_X$ this means that $Y \amalg_X Y$ exists. Since existence and description of pushouts is subtle in general, we have to impose some mild additional conditions.
\end{Rmk}

\begin{Prop}\label{Prop cogroup essentially surjective}
Let $X \in \Sch_S$, and let $(Y, c, \iota)$ be a cogroup in $X/\Sch_X$. Assume that the underlying set of $Y \amalg_X Y$ agrees with the pushout $|Y| \amalg_{|X|} |Y|$ of sets. Then $(Y,c,\iota) \cong (\NIL_X(\mathscr F),c,\iota)$ for a unique $\mathscr F \in \Qcoh(\mathcal O_X)$.
\end{Prop}

\begin{Rmk}\label{Rmk pushout separated}
The assumption on $Y \amalg_X Y$ is for example satisfied if $Y \to X$ is separated, or if formation of $Y \amalg_X Y$ commutes with base change of $X$. Indeed, in the former case, the section $\sigma \colon X \to Y$ is a closed immersion. Then the pushout $Y \amalg_X Y$ is described in \boldref{Rmk pushout of schemes} above, and agrees with the pushout as locally ringed spaces (in particular as sets).

In the latter case, it suffices to look at each fibre. But a section $\Spec k \to Y$ to a $k$-scheme is always closed, so we may proceed as in the first case. In fact, the first case follows from the second, because formation of $Y \amalg_X Y$ commutes with base change if $\sigma \colon X \to Y$ is closed, since the formation of $B \times_A B$ for a split surjection $B \twoheadrightarrow A$ commutes with tensor products.
\end{Rmk}

\begin{proof}[Proof of \boldref{Prop cogroup essentially surjective}.]
Because $Y \amalg_X Y$ is also the pushout as sets, we get a cogroup in $|X|/\Set_{|X|}$. Thus for each $x \in X$ we get a cogroup in $*/\Set$ (pointed sets). But these are always trivial: the existence of a two-sided counit $\varepsilon \colon Y \to *$ implies that the compositions
\begin{equation*}
\begin{tikzcd}
Y \ar{r}{c} & Y \vee Y \ar[yshift=.25em]{r}{1 \vee \varepsilon}\ar[yshift=-.25em]{r}[swap]{\varepsilon \vee 1} & Y
\end{tikzcd}
\end{equation*}
are the identity, where $\vee$ denotes the wedge sum. The equaliser of $1 \vee \varepsilon$ and $\varepsilon \vee 1$ is the point $*$, showing that $Y = *$. 
%
Applying this to all fibres $Y_x \to x$, we conclude that $Y \to X$ is a bijection, with section $\sigma$. Hence, $Y \to X$ is a homeomorphism, hence affine \cite[Tag \href{https://stacks.math.columbia.edu/tag/04DE}{04DE}]{Stacks}; say $Y = \SPEC_X(\mathscr A)$ for some quasi-coherent sheaf of $\mathcal O_X$-algebras $\mathscr A$. Then $\mathscr A$ is an abelian group object in $\Qcoh(\mathcal O_X)/\mathcal O_X$, say with multiplication $m = c^\#$ and retraction $\pi = \sigma^\#$. 

Write $\mathscr I = \ker(\mathscr A \stackrel\pi\to \mathcal O_X)$, so $\mathscr A \cong \mathcal O_X \oplus \mathscr I$ as $\mathcal O_X$-modules. We have natural maps
\[
\mathscr A \to \mathscr A \underset{\mathcal O_X}\times \mathscr A \cong \mathcal O_X \oplus \mathscr I \oplus \mathscr I
\vspace{-.25em}
\]
given by $(\pi,\id_\mathscr A)$ and $(\id_\mathscr A,\pi)$.
%
The relations of abelian groups imply that $m \circ (\id_\mathscr A,\pi) = \id_\mathscr A = m \circ (\pi,\id_{\mathscr A})$, so $m$ is necessarily given by
\begin{align*}
m \colon \mathcal O_X \oplus \mathscr I \oplus \mathscr I &\to \mathcal O_X \oplus \mathscr I\\
(a,b,c) &\mapsto (a,b+c).
\end{align*}
But $m$ is also an $\mathcal O_X$-algebra homomorphism. Thus,
\[
(0,f)\cdot(0,g) = m(0,f,0)\cdot m(0,0,g) = m((0,f,0)\cdot (0,0,g)) = 0,
\]
which means that $\mathscr I^2 = 0$. Thus, $(\mathscr A,m)$ is isomorphic to $(\Sym^{\leq 1} \mathscr I,m)$, where $m$ is the group structure induced by the addition $\mathscr I \oplus \mathscr I \to \mathscr I$. This shows that $(Y,c,\iota)$ is isomorphic to $(\NIL_X(\mathscr I),c,\iota)$. Since $\mathscr I = \ker(\mathscr A \to \mathcal O_X)$ can be recovered from $Y$, the uniqueness statement follows.
\end{proof}
\vskip-\lastskip
Write $\coAbsep(X/\Sch_X)$ for the category of abelian cogroup objects $(Y,c,\iota)$ in $X/\Sch_X$ such that the structure map $Y \to X$ is separated (see \boldref{Rmk pushout separated}).

\begin{Thm}\label{Thm Qcoh}
Let $X \in \Sch_S$. Then the functor
\begin{align*}
F \colon \Qcoh(\mathcal O_X)\op &\to \coAbsep(X/\Sch_X)\\
\mathscr F &\mapsto (\NIL_X(\mathscr F),c,\iota)
\end{align*}
is an equivalence of categories.
\end{Thm}

\begin{proof}
Essential surjectivity is \boldref{Prop cogroup essentially surjective}. If $\mathscr F, \mathscr G \in \Qcoh(\mathcal O_X)$, then an $\mathcal O_X$-algebra homomorphism $f \colon \mathcal O_X \oplus \mathscr F \to \mathcal O_X \oplus \mathscr G$ preserving the surjections to $\mathcal O_X$ necessarily maps $\mathscr F$ into $\mathscr G$, hence comes from a unique $\mathcal O_X$-module map $\mathscr F \to \mathscr G$. Moveover, any such $f$ automatically preserves the abelian cogroup structure, showing that $F$ is fully faithful.
\end{proof}

\section{Structure sheaf}\label{Sec structure sheaf}
We modify the argument of Brandenburg's second answer \cite{BraAns} to \cite{Bel} to categorically reconstruct the structure sheaf on $S$. The proof is based on the well-known formula
\[
\End(\id_{\Qcoh(\mathcal O_X)}) \cong \Gamma(X,\mathcal O_X).
\]
However, the left hand side does not naturally come with functorial restriction maps for morphisms $f \colon Y \to X$, for the same reason that sheaves $\mathscr F, \mathscr G$ on a topological space $X$ do not have restriction maps
\[
\Hom(\mathscr F(U),\mathscr G(U)) \to \Hom(\mathscr F(V),\mathscr G(V))
\]
for opens $V \subseteq U$. For sheaves, the solution is to work with $\Hom(\mathscr F|_U,\mathscr G|_U)$ instead of $\Hom(\mathscr F(U),\mathscr G(U))$, and our solution will be similar.

\begin{Def}\label{Def QCOH}
If $S$ is a scheme, write $\QCOH_{-/S}$ for the category whose objects are pairs $(X,\mathscr F)$, where $X$ is an $S$-scheme and $\mathscr F$ is a quasi-coherent sheaf on $X$.
A morphism $(f,\phi) \colon (X,\mathscr F) \to (Y,\mathscr G)$ consists of a morphism $f \colon X \to Y$ of $S$-schemes and a morphism $\phi \colon f^*\mathscr G \to  \mathscr F$ of $\mathcal O_X$-modules (equivalently, a morphism $\mathscr G \to f_*\mathscr F$ of $\mathcal O_Y$-modules). 
The forgetful functor $p \colon \QCOH_{-/S} \to \Sch_S$ makes it into a fibred category \cite[Tag \href{https://stacks.math.columbia.edu/tag/03YM}{03YM}]{Stacks}. The fibre $\QCOH_{X/S} = p^{-1}(\id_X)$ is canonically equivalent to $\Qcoh(\mathcal O_X)\op$.

Similarly, define the category $\QCOH_{-/-/S}$ whose objects are pairs $(Y \to X, \mathscr F)$, where $Y \to X$ is a morphism of $S$-schemes and $\mathscr F$ is a quasi-coherent sheaf on $Y$. A morphism $(f,g,\phi) \colon (Y \to X,\mathscr F) \to (Y' \to X', \mathscr F')$ in $\QCOH_{-/-/S}$ consists of morphisms $f \colon Y \to Y'$ and $g \colon X \to X'$ of $S$-schemes such that the diagram
\begin{equation}\label{Dia morphism in QCOH}
\begin{tikzcd}
Y \ar{r}{f}\ar{d} & Y'\ar{d}\\
X \ar{r}{g} & X'
\end{tikzcd}
\end{equation}
commutes, along with a morphism $\phi \colon f^*\mathscr F' \to \mathscr F$ of $\mathcal O_Y$-modules. It has a forgetful functor $q \colon \QCOH_{-/-/S} \to \Sch_S$ mapping $(Y \to X,\mathscr F)$ to $X$, making it into a fibred category whose fibre $\QCOH_{-/X/S} = q^{-1}(\id_X)$ is canonically equivalent to $\QCOH_{-/X}$. For a morphism $f \colon Y \to X$ in $\Sch_S$, define the pushforward functor
\begin{align*}
f_* \colon \QCOH_{-/Y/S} &\to \QCOH_{-/X/S}\\
(Z \to Y, \mathscr F) &\mapsto (Z \to X, \mathscr F),\\
(g,\id_Y,\phi) &\mapsto (g,\id_X,\phi).
\end{align*}
The pushforward $f_* \colon \QCOH_{-/Y/S} \to \QCOH_{-/X/S}$ allows us to turn the association $X \mapsto \End(\id_{\QCOH_{-/X/S}})$ into a functor $\End(\id)$: for a morphism $f \colon Y \to X$ in $\Sch_S$, define the pullback
\[
f^* \colon \End(\id_{\QCOH_{-/X/S}}) \to \End(\id_{\QCOH_{-/Y/S}})
\]
by $f^*(\alpha)_y = \alpha_{f_*y}$ for any $\alpha \in \mathscr \End(\id_{\QCOH_{-/X/S}})$ and any $y \in \QCOH_{-/Y/S}$.
\end{Def}

\begin{Def}\label{Def E}
Let $\mathscr E \subseteq \End(\id)$ be the subfunctor consisting of $\alpha \in \End(\id)$ such that $\alpha_{(Z \to X,\mathscr F)} \colon (Z \to X,\mathscr F) \to (Z \to X,\mathscr F)$ is of the form\footnote{A general $\alpha_{(Z \to X, \mathscr F)}$ has the form $(f,\id_X,\phi)$, so the extra assumption is $f = \id_Z$.} $(\id_Z,\id_X,\phi)$, for all $(Z \to X,\mathscr F)$ in $\QCOH_{-/-/S}$. This property is clearly preserved by the pullback $f^*$, so it defines a subfunctor. 

Concretely, an element $\alpha \in \mathscr E(X)$ consists of the following data: for every object $(Z \to X,\mathscr F) \in \QCOH_{-/X/S}$, a map
\begin{equation}
\alpha_{(Z \to X,\mathscr F)} = (\id_Z,\id_X,\phi_{(Z \to X,\mathscr F)}) \colon (Z \to X,\mathscr F) \to (Z \to X,\mathscr F)\label{Eq pullback alpha}
\end{equation}
in $\QCOH_{-/X/S}$ for some morphism of $\mathcal O_Z$-modules $\phi_{(Z \to X,\mathscr F)} \colon \mathscr F \to \mathscr F$, such that for every map $(f,\id_X,\psi) \colon (Z \to X,\mathscr F) \to (Z' \to X,\mathscr F')$ in $\QCOH_{-/X/S}$, the diagram
\begin{equation}\label{Dia End}
\begin{tikzcd}[column sep=5em]
f^* \mathscr F' \ar{r}{\phi_{(Z' \to X,\mathscr F')}}\ar{d}[swap]{\psi} & f^* \mathscr F' \ar{d}{\psi}\\
\mathscr F \ar{r}{\phi_{(Z \to X,\mathscr F)}} & \mathscr F
\end{tikzcd}
\end{equation}
commutes. 
%
If $f \colon Y \to X$ is a morphism, then the pullback $f^* \colon \mathscr E(X) \to \mathscr E(Y)$ is defined by sending $\alpha \in \mathscr E(X)$ to $f^* \alpha \colon \id_{\QCOH_{-/Y/S}} \to \id_{\QCOH_{-/Y/S}}$ given by
\[
\left(f^*\alpha\right)_{(Z \to Y,\mathscr F)} = (\id_Z,\id_Y,\phi_{(Z \to X,\mathscr F)}) \colon (Z \to Y,\mathscr F) \to (Z \to Y,\mathscr F)
\]
for all $(Z \to Y,\mathscr F) \in \QCOH_{-/Y/S}$. Thus, in (\ref{Eq pullback alpha}) we may unambiguously write $\phi_{\mathscr F}$ instead of $\phi_{(Z \to X,\mathscr F)}$, because it does not depend on which $X$ we consider (and we think of $Z$ as understood). We think of $\mathscr E(X)$ as the association of an endomorphism $\phi_{\mathscr F}$ for every quasi-coherent sheaf $\mathscr F$ on an $X$-scheme $Z$, with the compatibility condition (\ref{Dia End}).
\end{Def}

\begin{Prop}\label{Prop structure sheaf}
Let $S$ be a scheme. Then the presheaf $\mathscr E \colon \Sch_S\op \to \Ring$ is isomorphic to $\mathcal O \colon X \mapsto \Gamma(X,\mathcal O_X)$.
\end{Prop}

\begin{proof}
We will define maps $a \colon \mathcal O \to \mathscr E$ and $b \colon \mathscr E \to \mathcal O$ that are mutual inverses. To define $a$, let $X \in \Sch_S$ be given. If $f \in \mathcal O(X)$, then define $a(f) \in \mathscr E(X)$ by letting $\phi_{\mathscr F}$ for a sheaf $\mathscr F$ on $g \colon Z \to X$ be multiplication by $g^\# f \in \mathcal O(Z)$. Clearly the diagrams (\ref{Dia End}) commute for every morphism in $\QCOH_{-/X/S}$, so $a(f)$ defines an element of $\mathcal E(X)$. Similarly, the definition of the pullback $f^* \colon \mathscr E(X) \to \mathscr E(Y)$ immediately implies that $a$ is a morphism of presheaves.

Conversely, if $X \in \Sch_S$ and $\alpha \in \mathscr E(X)$, then set $b(\alpha) = \phi_{\mathcal O_X}(1) \in \mathcal O_X(X)$. For any $f \colon Y \to X$, commutativity of (\ref{Dia End}) for $\alpha \in \mathscr E(X)$ and for the morphism $(f,f^\#) \colon (Y,\mathcal O_Y) \to (X,\mathcal O_X)$ gives
\[
f^*\phi_{\mathcal O_Y}(1) = \phi_{\mathcal O_Y}(1) = f^\#(b(\alpha)),
\]
showing that $b$ is a morphism of presheaves.

It is clear that $ba(f) = f$ for any $f \in \mathcal O(X)$. Conversely, given $X \in \Sch_S$ and $\alpha \in \mathscr E(X)$, we must show that $\phi_{\mathscr F}$ is multiplication by $f = b(\alpha)$ for any quasi-coherent sheaf $\mathscr F$ on an $X$-scheme $Z$.
%
Sections $\mathscr F(U)$ for $U \subseteq Z$ open are given by maps $\psi \colon \mathcal O_U \to \mathscr F|_U$, so the diagram (\ref{Dia End}) for the morphisms $(U \to X,\id_X,\psi) \colon (U \to X,\mathscr F) \to (X \to X,\mathcal O_X)$ shows that each element of $\phi_{\mathscr F}(U)$ gets multiplied by $f$.
\end{proof}

\begin{Thm}\label{Thm forgetful functor}
Let $S$ be a scheme. Then the forgetful functor $U \colon \Sch_S \to \Sch$ can be reconstructed categorically from $\Sch_S$.
\end{Thm}

\begin{proof}
By \boldref{Prop topology}, the topology can be reconstructed categorically. Moreover, the fibred category $\QCOH_{-/-/S} \to \Sch_S$ can be reconstructed categorically: by \boldref{Thm Qcoh} it is equivalent to the category of pairs $(Y \to X, \mathscr F)$ where $Y \to X$ is a morphism in $\Sch_S$, and $\mathscr F$ a separated abelian cogroup object of $Y/\Sch_Y$. Then the presheaf $\mathscr E \colon \Sch_S\op \to \Ring$ of \boldref{Def E} can be reconstructed categorically, which by \boldref{Prop structure sheaf} is isomorphic to the sheaf $\mathcal O$ on the big Zariski site $\Sch_S$. 

This shows that the functor $\Sch_S \to \LRS$ given by $X \mapsto (|X|,\mathcal O_X)$ can be reconstructed categorically. Its essential image lands in the category of schemes.
\end{proof}

\section{Proofs of main theorems}\label{Sec proofs of main theorems}
With the results from \boldref{Sec structure sheaf}, we are ready to prove \boldref{Thm main} and \boldref{Thm rings}.

\begin{proof}[Proof of \boldref{Thm main}]
By \boldref{Thm forgetful functor}, the forgetful functor $U \colon \Sch_S \to \Sch$ can be reconstructed categorically from $\Sch_S$. Then \boldref{Lem categorical reconstruction equivalence on ISOM} implies that the functor
\begin{align*}
\Isom(S,S') &\to \ISOM(\Sch_{S'},\Sch_S)\\
f &\mapsto f^*
\end{align*}
is an equivalence.
\end{proof}

\begin{proof}[Proof of \boldref{Thm rings}]
Analogously to the construction of $\QCOH_{-/-/S}$, we can form a category $\Mod_{-/-/R}$ of modules over $R \to A \to B$, and use the pushforward
\[
f_* \colon \Mod_{-/B/R} \to \Mod_{-/A/R}
\]
to make $A \mapsto \End(\id_{\Mod_{-/A/R}})$ functorial in $A$. Similarly, we construct a subfunctor
\[
\mathscr E \colon \Alg_R \to \Ring
\]
of endomorphisms $\alpha \colon \id_{\Mod_{-/A/R}} \to \id_{\Mod_{-/A/R}}$ that are the identity on the rings $B/A/R$ (so they only act on the module). Analogously to \boldref{Prop structure sheaf}, one shows that $\mathscr E$ is isomorphic to the forgetful functor $\Alg_R \to \Ring$. Then \boldref{Lem categorical reconstruction equivalence on ISOM} gives the result.
\end{proof}
\vskip-\lastskip
We note that the proof of \boldref{Thm rings} only relies on the constructions from \boldref{Sec structure sheaf}, so it is much shorter than the proof of Clark--Bergman \cite{CB}, who only treat the case where the base is an integral domain. Already in this case, the statement that $\ISOM(\Alg_R,\Alg_{R'})$ is a setoid appears to be new.

\phantomsection
\printbibliography
\end{document}